\documentclass[12pt]{article}
\usepackage{amsmath,amsfonts,graphicx}
\usepackage[framed,numbered,autolinebreaks,useliterate]{mcode}
\usepackage{tabularx}
\usepackage{geometry} 
\newcommand{\argmin}[1]{\underset{#1}{\mathrm{argmin~}}}

\newcommand{\Rn}{\mathcal{R}_P}
\newcommand{\rank}{\mathrm{rank}}
\newcommand{\real}{\mathsf{Re}}
\newcommand{\realR}{\mathbb{R}}
\newcommand{\imag}{\mathsf{Im}}
\newcommand{\imagI}{\mathbb{I}}
\newcommand{\ii}{\mathrm{i}}

\newcommand{\C}{\mathbb{C}}

\newcommand{\freq}{\nu_p}
\newcommand{\damp}{\gamma_p}

\graphicspath{{Figures/}}

\title{Frequency estimation based on Hankel matrices and the alternating direction method of multipliers}
\usepackage{authblk}

\title{Frequency estimation based on Hankel matrices and the alternating direction method of multipliers}
\author[1]{Fredrik Andersson\thanks{fa@maths.lth.se}}
\author[1]{Marcus Carlsson\thanks{mc@maths.lth.se}}
\author[2]{Jean-Yves Tourneret\thanks{jean-yves.tourneret@enseeiht.fr}}
\author[2]{Herwig Wendt\thanks{herwig.wendt@irit.fr}}
\affil[1]{\small Centre for Mathematical Sciences, Lund University, Lund, Sweden}
\affil[2]{\small University of Toulouse, CNRS UMR 5505, INP-ENSEEIHT, Toulouse, France}

\begin{document}

\maketitle
\begin{abstract}
We develop a parametric high-resolution method for the estimation of the frequency nodes of linear combinations of complex exponentials with exponential damping. We use Kronecker's theorem to formulate the associated nonlinear least squares problem as an optimization problem in the space of vectors generating Hankel matrices of fixed rank. Approximate solutions to this problem are obtained by using the alternating direction method of multipliers. Finally, we extract the frequency estimates from the con-eigenvectors of the solution Hankel matrix. The resulting algorithm is simple, easy to implement and can be applied to data with equally spaced samples with approximation weights, which for instance allows cases of missing data samples. By means of numerical simulations, we analyze and illustrate the excellent performance of the method, attaining the Cram\'er-Rao bound.
\end{abstract}


%
%

\section{Introduction}
\label{sec:intro}
Spectral estimation constitutes a classical problem that has found applications in a large variety of fields (including astronomy, radar, communications, economics, medical imaging, spectroscopy, \dots, to name but a few).
One important category of spectral estimation problems arises for  signals that can be well represented by the parametric model
\begin{equation}
\label{equ:sexpo}
f_c(t)= \sum_{p=0}^{P-1} c_p e^{\zeta_p t} + e(t),\quad c_p,\zeta_p\in \C,
\end{equation}
where $e(t)$ is an additive noise term.  Given a vector of (typically equally spaced) samples,
$$
f(j)=f_c\left(t_0+j t_s\right),\quad0 \le j \le 2 N
$$
where $t_s$ is the sampling period, the goal is to estimate  the complex frequency nodes
$
\zeta_p=2\pi(\damp+\mbox{i}\freq),
$
i.e., the damping and frequency parameters $\damp$  and  $\freq$.
Note that once the parameters $\boldsymbol{\zeta}=(\zeta_0,\cdots,\zeta_{P-1})^T$ have been computed, determining $\mathbf{c}=(c_0,\cdots,c_{P-1})^T$ reduces to a simple linear regression problem. Thus, the focus lies on the estimation of the nodes $\boldsymbol{\zeta}$.
The signal model \eqref{equ:sexpo} includes the case of exponentially damped signals defined by $\damp<0$. This case has recently received significant interest, notably in applications involving nuclear quadrupole resonance signals for which the $\damp$ (together with $\freq$) are used as a signature of the chemical composition of an analysed substance (see for instance \cite{Gudmundson2012} and references therein).

The literature on methods for estimating  $\boldsymbol{\zeta}$ in this setting is rich, cf. \cite{Stoica2005} for an overview. Here, we focus on so-called \emph{high-resolution} methods, i.e., methods whose frequency resolution is not tied to a pre-defined grid (and can attain numerical precision in the noise free case).
\\
One important class of techniques is based on the principle of nonlinear least squares (NLS) \cite{NLLS1986} and aims at solving the nonlinear regression problem associated with \eqref{equ:sexpo}, i.e., at minimizing $$\theta(\boldsymbol{c},\boldsymbol{\zeta})=\sum_{j=0}^{2N}\big| f(j)-  \sum_{p=0}^{P-1} c_p e^{\zeta_p j} \big|^2$$ with respect to the parameter vectors $\boldsymbol{c}$ and $\boldsymbol{\zeta}$. Although NLS enjoy desirable theoretical properties (e.g., those of maximum likelihood estimators when $e(t)$ is a Gaussian white noise; robustness to coloured noise), their practical use remains limited due to the extreme difficulty of globally minimizing $\theta(\boldsymbol{c},\boldsymbol{\zeta})$ due to its pronounced multimodal shape with one very sharp global minimum  (cf., e.g.,  \cite{Stoica2005} and references therein).
\\
A second prominent and popular class is given by subspace methods (such as  MUSIC \cite{Schmidt}, ESPRIT \cite{ESPRIT} and min-norm \cite{MinNorm1983}). These methods are based on estimates of the covariance of $f$ and rely on the assumption that the noise $e$ is white.  

In the present contribution, we develop a novel high-resolution methodology for the estimation of $\boldsymbol{\zeta}$ in \eqref{equ:sexpo}. Similar to NLS, we aim at approximating $f$ as good as possible by a linear combination of complex exponentials. However, the proposed methodology fundamentally departs from any of the above methods in the following ways:

First, it is based on \emph{Kronecker's theorem} for complex symmetric matrices, which essentially states that if a function $f_c$ is uniformly sampled, then the \emph{Hankel matrix} that is generated by the vector of samples $f$ has rank $P$ if and only if $f_c$ coincides at the sample points with a function that is a linear combination of $P$ exponential functions. This fact has been used for the sparse approximation of functions by sums of exponentials in the con-eigenvalue approach \cite{Beylkin_Monzon_2005} and the alternating projections method \cite{AlternatingProjections}. Here, it is used to formulate an NLS type minimization problem in which the model (and its parameters $\boldsymbol{c}$ and $\boldsymbol{\zeta}$) does not enter explicitly. Instead, the model is imposed \emph{implicitly} by constraining the rank of the Hankel matrix generated by the vector $g$ approximating $f$ to equal $P$. Consequently, residual minimization is not performed  over the parameter space $\{\boldsymbol{c},\boldsymbol{\zeta}\}$ directly, but over the space of vectors $g$ which generate Hankel matrices of rank $P$. The frequency nodes $\boldsymbol{\zeta}$ are then obtained by considering the con-eigenvectors of the solution Hankel matrix.

Second, we reformulate the minimization problem such that it can be effectively solved by the \emph{alternating direction method of multipliers} (ADMM) \cite{Boyd_ADMM}. ADMM is an iterative technique that is recently gaining popularity due to its robustness, versatility and applicability to large-scale problems. Moreover, this technique enjoys performance comparable to problem-specific state-of-the-art methods. While the optimization problem considered here is nonconvex, it is shown numerically that the solutions obtained with the ADMM procedure generically provide excellent approximations and parameter estimates for the model \eqref{equ:sexpo}.

The resulting Hankel matrix ADMM frequency estimation procedure is practically extremely appealing due to its simplicity and ease of implementation. It enjoys excellent performance, outperforming subspace methods while at the same time alleviating the practical limitations of classical NLS.
Unlike subspace methods, it does not rely on a specific noise model.
Another interesting property is that it applies to situations with \emph{missing} (or censored) data.

The remainder of this work is organized as follows.
In Section \ref{sec:admm}, we define an optimization problem based on Kronecker's theorem for the approximation by sums of complex exponentials and formulate the ADMM-based procedure for obtaining its solution.
Section \ref{sec:freqest} summarizes the procedure for obtaining frequency node estimates from the solution Hankel matrix.
In Section \ref{sec:results}, we analyze the performance of the method by means of numerical simulations.
Section \ref{sec:conclusions} concludes this contribution and points to future work.

\section{Approximation by sums of exponentials using ADMM}
\label{sec:admm}

In this section, we approximate the vector $f$ as good as possible by a vector $g$ that is a linear combination of complex exponential functions. Using Kronecker's theorem, we express the solution to this approximation problem in terms of the vector $g$ whose Hankel matrix has rank $P$ and minimizes the residual. The resulting optimization problem will then be reformulated and solved by ADMM.

\subsection{\bf Hankel matrices and Kronecker's theorem}
For completeness, we recall that a {Hankel} matrix $A\in\mathcal{H}$%
 has constant values on the anti-diagonals, i.e., $$A(j,k)=A(j',k') \mbox{ if $j+k=j'+k'$}.$$
It can thus be generated elementwise from a vector $g(j)$, such that 
$$A(j,k)=g(j+k), \quad 0\le j,k\le N.$$ 
We denote this by $A=Hg$.

\emph{Kronecker's theorem} \cite{AAK, JAT} states that if the Hankel matrix $A=Hg$ generated by the vector $g$ is of rank $P$ then, with the exception of degenerate cases, there exists $\{\zeta_p\}_{p=0}^{P-1}$ and $\{c_p\}_{p=0}^{P-1}$ such that
\begin{equation}\label{Kronecker}
g(j) = \sum_{p=0}^{P-1} c_p e^{\zeta_p j}, \quad c_p,\zeta_p\in \C.
\end{equation}
The converse holds as well.

It follows that the best approximation (in the $l_2$ sense) of $f$ by a linear combination of $P$ complex exponentials is given by the vector $g$ which satisfies $\rank(Hg)= P$ and minimizes the $l_2$ norm of the residual $r=f-g$. In other words, $g$ can be obtained as the solution of the optimisation problem
\begin{equation} \label{opt0}
\begin{aligned}
& \underset{g}{\text{minimize}}
& & \frac{1}{2}\| r\|_2^2 = \frac{1}{2} \| f-g\|_2^2 \\
& \text{subject to}
& & \rank(Hg)= P.
\end{aligned}
\end{equation}

The parameter $P$ can be selected by considering the singular values of the Hankel matrix $Hf$ generated by the samples $f$, cf., \cite{Beylkin_Monzon_2005}.
At this stage, we assume $P$ to be given.

\subsection{A solution based on ADMM}

Let $\Rn$ be the  indicator function for square matrices $S$ given by
$$
\Rn(S) = \begin{cases}
0 & \mbox{if } \rank(S) = P, \\
\infty & \mbox{otherwise.}
\end{cases}
$$
Then, problem \eqref{opt0} can be reformulated as follows
\begin{equation} \label{hankel_eq}
\begin{aligned}
& \underset{A,r}{\text{minimize}}
& & \Rn(A) + \frac{1}{2} \| r\|_2^2 \\
& \text{subject to}
& & A(j,k) + r(j+k)= f(j+k), \quad 0 \le j,k \le N.
\end{aligned}
\end{equation}
The optimization problem \eqref{hankel_eq} has a cost function that consists of two terms, each depending only on one of the variables $A$ and $r$, along with a linear constraint. The problem formulation is therefore well suited to be addressed using ADMM. ADMM is an iterative technique in which a solution to a large global problem is found by coordinating solutions to smaller subproblems. 
It can be seen as an attempt to merge the benefits of dual decomposition and augmented Lagrangian methods for constrained optimization.
For an overview of the ADMM method see \cite{Boyd_ADMM}. For a convex cost function ADMM is guaranteed  to converge to the unique solution. For non-convex problems it can get stuck at local minima, although it can work quite well in practice also in these situations, cf. \cite[Chapter 9]{Boyd_ADMM}. Since the rank constraint $\Rn$ is nonconvex, convergence of an algorithm solving \eqref{hankel_eq} is not guaranteed. Our numerical simulations indicate that the method works substantially better than established high resolution techniques like ESPRIT, and moreover that it can be applied to cases where ESPRIT is non-applicable (e.g., missing data)\footnote{A convex problem could be obtained by replacing the $\Rn$ by the nuclear norm, at the cost of biased solutions.}.

An ADMM iterative step (from iteration $q$ to iteration $q+1$) for (\ref{hankel_eq})  reads
\begin{align}
\label{admm1}& A^{q+1} = \argmin{A} L(A,r^q,\Lambda^q),\\
\label{admm2}&  r^{q+1} = \argmin{r} L(A^{q+1},r,\Lambda^q),\\
& \Lambda^{q+1}(j,k) = \Lambda^q(j,k) + \rho \big( A^{q+1}(j,k) \nonumber\\
& \qquad\qquad\qquad\quad\qquad+r^{q+1}(j+k)-f(j+k) \big)
\end{align}
where $\Lambda = \Lambda_\realR + \ii \Lambda_\imagI$ are the Lagrange multipliers,
and $L(A,r,\Lambda)$ is the augmented Lagrangian associated with (\ref{hankel_eq}),
\begin{multline*}
L(A,r,\Lambda) = \Rn(A) + \frac{1}{2} \| r\|_2^2 + \\
\sum_{j,k=0}^{N} \Big[  \Lambda_\realR(j,k) \real\big(A(j,k)+r(j+k)-f(j+k) \big) \\
\;\;\;\quad+ \Lambda_\imagI(j,k) \imag\big(A(j,k)+r(j+k)-f(j+k) \big)
 \\+ \frac{\rho}{2} \big |A(j,k)+r(j+k)-f(j+k) \big |^2 \Big].
\end{multline*}
The first minimization step \eqref{admm1} is nonconvex due to the projection onto a non convex set but can be solved analytically. Deriving closed form expressions for both minimization steps steps \eqref{admm1} and \eqref{admm2} is straightforward.
The solution $A^{q+1}$ to the first minimization step is equal to the best rank $P$ approximation of the matrix $B$ defined elementwise by
$$
B(j,k) = f(j+k)-r^q(j+k)-\rho^{-1} \Lambda^q(j,k).
$$
Denoting by $B=U \Sigma V^\ast$ the singular value decomposition of $B$, then by the Eckart-Young theorem,
\begin{equation}
A^{q+1} = U \Sigma_P V^\ast,
\end{equation}
where $\Sigma_P$ is obtained from the diagonal matrix $\Sigma$ by replacing the diagonal elements $\Sigma(l,l)$ for $l>P$ with zeros.
\\
The solution to the second minimization step is given by
\begin{equation}
\label{equ:step2}
r^{q+1}(l) = \frac{\rho Q(l) f(l)-\sum_{j+k=l} \big[\Lambda^q(j,k) +{\rho} A^{q+1}(j,k) \big]}{1+ \rho Q(l) },
\end{equation}
where $\sum_{j+k=l} $ stands for $\sum_{\{0\le j,k \le N\} \cap \{j+k=l\}} $ and
$$
Q(l) = \begin{cases}
l+1 & \mbox{if } l \le N,\\
2N+1-l & \mbox{otherwise}.
\end{cases}
$$
With these explicit expressions for the ADMM minimization steps, the procedure for the approximation by sums of exponentials is given by the following simple MATLAB function:
\vskip-0.75\baselineskip
\begin{lstlisting}
function f_appro=expo_sum_admm(f,P,rho,it_max)
 %Best approximation to f using P exponentials.
 N=(length(f)-1)/2;r=zeros(2*N+1,1);
 Lambda=zeros(N+1);Q=[1:N+1,N:-1:1]';
 for iter=1:it_max,
  B=form_hank(f-r)-Lambda/rho;[u,s,v]=svd(B);
   s=diag(s);s(P+1:end)=0;A=u*diag(s)*v';
  r=(rho*Q.*f-sum_hank(Lambda+rho*A))./(1+rho*Q);
  Lambda=Lambda+rho*(A-form_hank(f-r));
 end;
 f_appro=sum_hank(A)./Q;
 
function H=form_hank(f)     %form Hankel matrix
 H=hankel(f(1:(end+1)/2),f((end+1)/2:end));
function f=sum_hank(H)      %sum anti-diagonals
 N=size(H,1)-1; f=zeros(2*N+1,1); H=flipud(H);
 for j=-N:N, f(j+N+1,1)=sum(diag(H,j)); end;
\end{lstlisting}

\subsection{Extension to the missing data case}
It is straightforward to replace $||r||_2^2$  in problem \eqref{hankel_eq} by a weighted norm
$$
||r||_{2,w}^2=\sum_{j=0}^{2N}w(j)|r(j)|^2
$$
and to derive the corresponding ADMM procedure. The only change in the above final expressions is that the denominator in the second minimization step \eqref{equ:step2} (and in the corresponding line 8 of the MATLAB code) is replaced by $w(l)+\rho Q(l)$.

Let $\mathcal{J}$ denote the set of indices $j$ of the missing data samples and set $f(\mathcal{J})=0$. The use of weights $w_\mathcal{J}$ defined by
$$
w(j)=w_\mathcal{J}(j)=\begin{cases}0&\mbox{if }j\in\mathcal{J}\\1&\mbox{otherwise}\end{cases}
$$
yields an ADMM procedure for data with missing samples. 



\medskip

\subsection{Convergence}
In order to highlight the simplicity of the proposed algorithm, the ADMM procedure was described above to terminate after a preset  number $q^*$ of iterations. Commonly used stopping criteria are based on the magnitude of the primal and dual residuals and are straight-forward to incorporate in the ADMM procedure (cf.,  \cite{Boyd_ADMM} for details).

The ADMM procedure would converge to the global solution of the problem \eqref{hankel_eq} if  the problem was convex. Since it is non convex, it is only guaranteed to converge locally and can converge to different points depending on the choice of $\rho$ and initial points $\Lambda^0$ and $r^0$.  Numerical experiments indicate, however, that the ADMM solution $g=f-r^{q^*}$ of the problem \eqref{hankel_eq} provides in general an excellent approximation for $f$. Here, we have initialized $\Lambda^0$ and $r^0$ with zeros.


\section{Frequency estimation}
\label{sec:freqest}
Once the solution Hankel matrix $A=Hg$ to \eqref{hankel_eq} has been found, the parameter vector $\boldsymbol{\zeta}$ can be obtained by considering the so-called {con-eigenvectors} of $A$.
Hankel matrices belong to the class of complex symmetric matrices satisfying $A=A^T$, which generically can be decomposed as \cite{Horn}
\begin{equation}
A = \sum_{p=0}^{N} s_p u_p u_p^T, \quad s_p\in \realR^+, \quad u_p \in \C^{N+1},
\end{equation}
where $s_0\geq s_1\geq\cdots\geq s_{N}>0$ are the decreasingly ordered con-eigenvalues of $A$ and where  the \emph{con-eigenvectors} $\{u_p\}$ are orthogonal and satisfy the relation
$A \overline{u_p} = s_p u_p$.
Similarly to the Eckart Young theorem, it can be shown that the best rank $P\leq N+1$ approximation of $A$ is given by \cite{Horn}
\begin{equation} \label{rank_n_rep}
\sum_{p=0}^{P-1} s_p u_p u_p^T.
\end{equation}

Since $g$ is a sum of $P$ exponentials, $Hg$ can be expressed in the form (\ref{rank_n_rep}), and therefore each con-eigenvector $u_p$ will also be a sum of the same exponentials. Let $U = (u_0, \dots, u_{P-1})$. It is then possible to write $U=V G$, where $V$ is the $N+1\times P$ Vandermonde matrix generated by $e^{\zeta_p}$ (i.e., $V(j,p)=e^{\zeta_p j}$) and $G$ is some (invertible) $P\times P$ matrix.
Denote as $\underline{U}$ (resp. $\overline{U}$) the matrix $U$ whose first row (resp. last row) has been dropped. Clearly, we have $\underline{U} = \underline{V} G$, $\overline{U} = \overline{V} G$ and the Vandermonde structure of $V$ leads to $\mathrm{diag}( e^{\zeta_0},\dots,e^{\zeta_{P-1}}) \overline{V} = \underline{V}.$
It follows that
$$
(\overline{U}^\ast \overline{U})^{-1} (\overline{U}^\ast \underline{U}) = \overline{U}^\dagger \underline{U} = G^{-1} \mathrm{diag}( e^{\zeta_0},\dots,e^{\zeta_{P-1}}) G,
$$
where $^\ast$ stands for conjugate transpose, where $~^\dagger$ denotes the pseudo-inverse operator.
Therefore, we can compute the nodes $\boldsymbol{\zeta}$ by computing the eigenvalues of $\overline{U}^\dagger \underline{U}$. 

\setlength{\tabcolsep}{1.0mm}
\begin{table}
\centering
\renewcommand{\arraystretch}{0.9}
\begin{tabular}{|r||r|r|r|r|r|r|r|}
\hline
$\freq$&$1.86$&$6.59$&$7.49$&$19.84$&&&\\
$\damp$&$0.40$&$-0.56$&$0.08$&$-0.45$&&&\\
$|c_p|$&$1.00$&$0.40$&$1.50$&$0.70$&&&\\\hline\hline
$\freq$&$1.86$&$    3.84$&$    5.95$&$    7.49$&$    9.60$&$   19.84$&$   30.08$\\
$\damp$&$   0.05$&$  0.14$&$   0.06$&$   0.01$&$   0.16$&$   0.08$&$  0.07$\\
$|c_p|$&$1.00$&$0.40$&$1.50$&$0.70$&$0.60$&$1.20$&$1.00$\\\hline
\end{tabular}
\caption{\label{tab:param}Parameters for the $P=4$ (top) and $P=7$ (bottom) complex exponentials used in the numerical simulations.}
\end{table}
\begin{figure}[tb]
\centering
\includegraphics[width=0.5\linewidth]{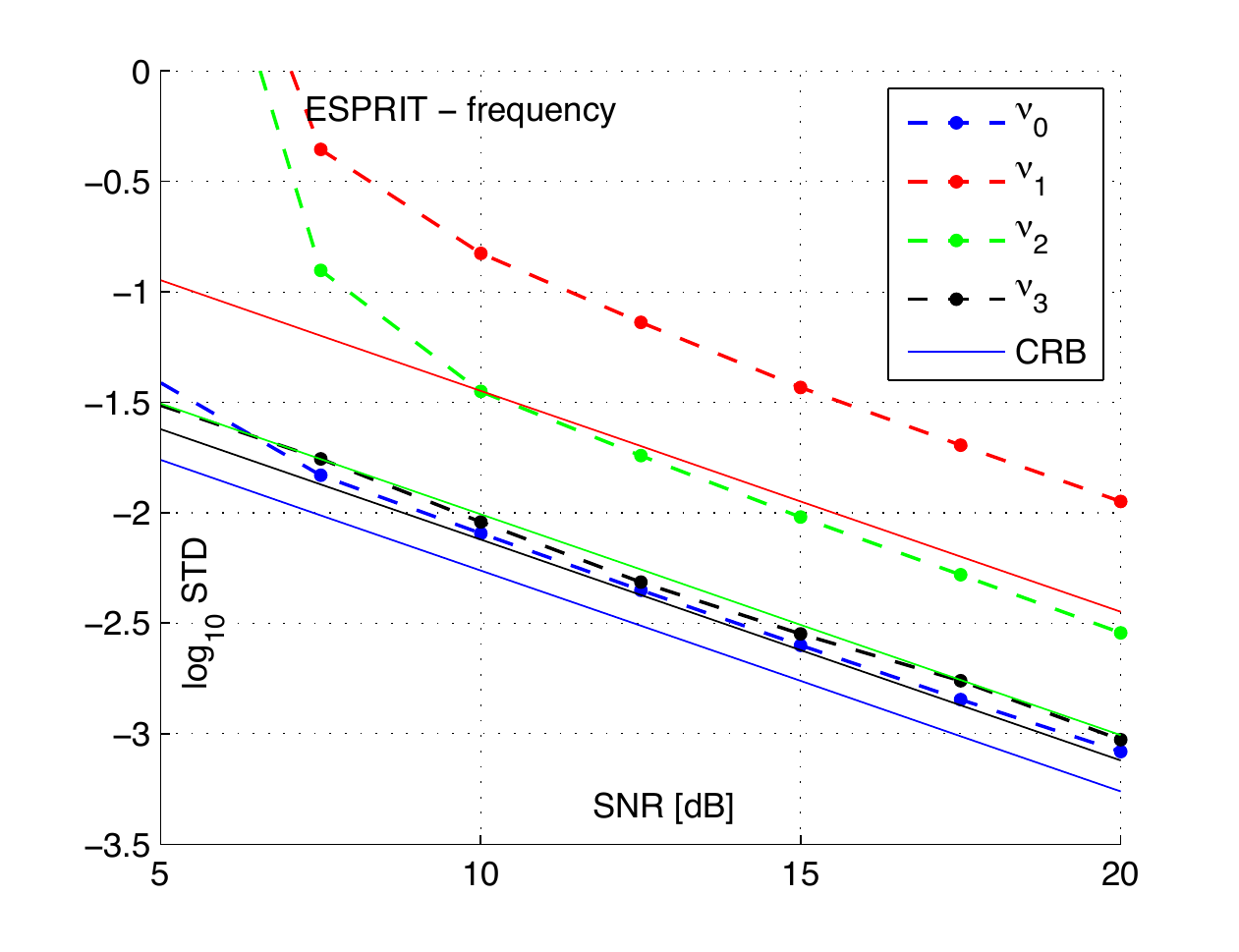}\includegraphics[width=0.5\linewidth]{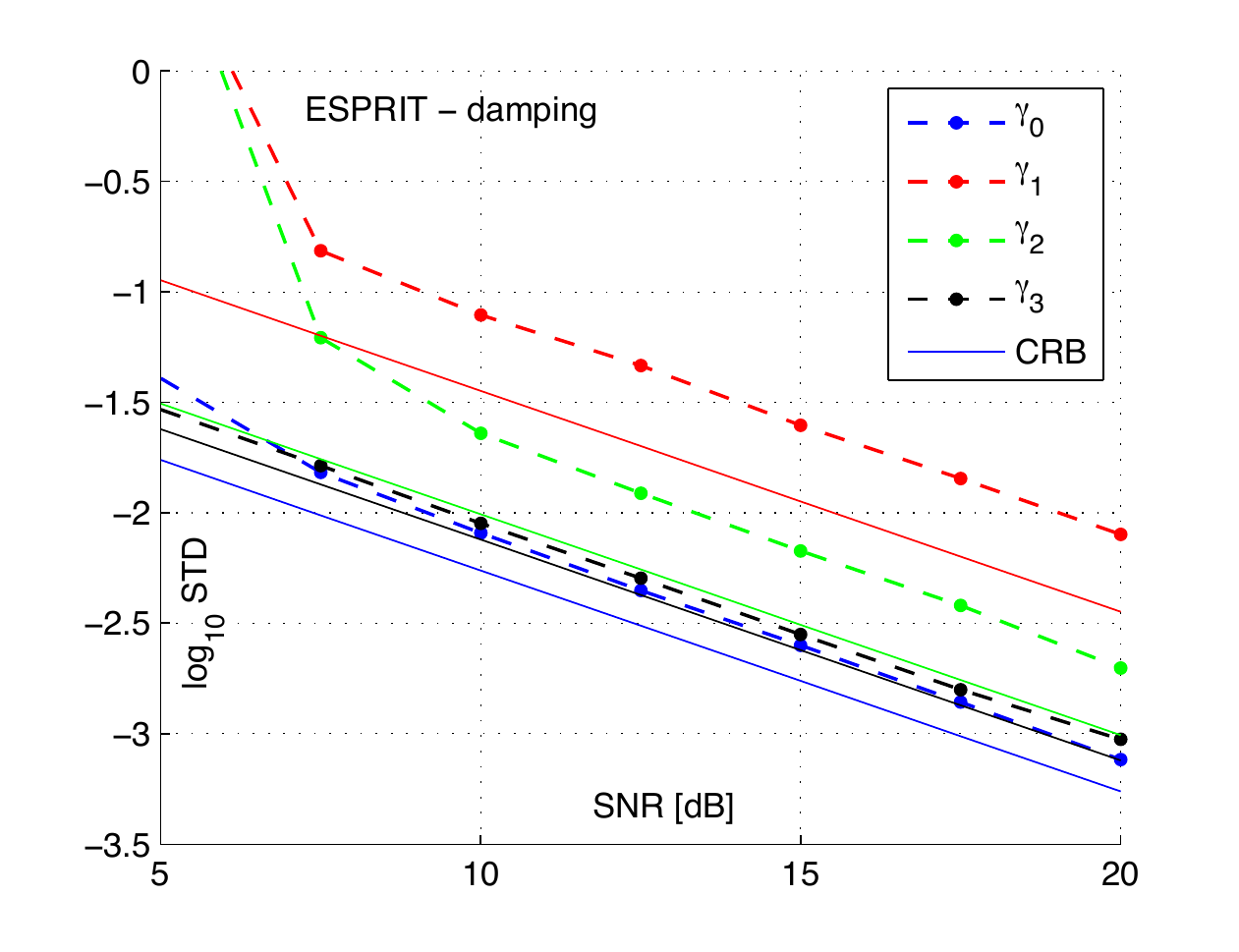}\\
\includegraphics[width=0.5\linewidth]{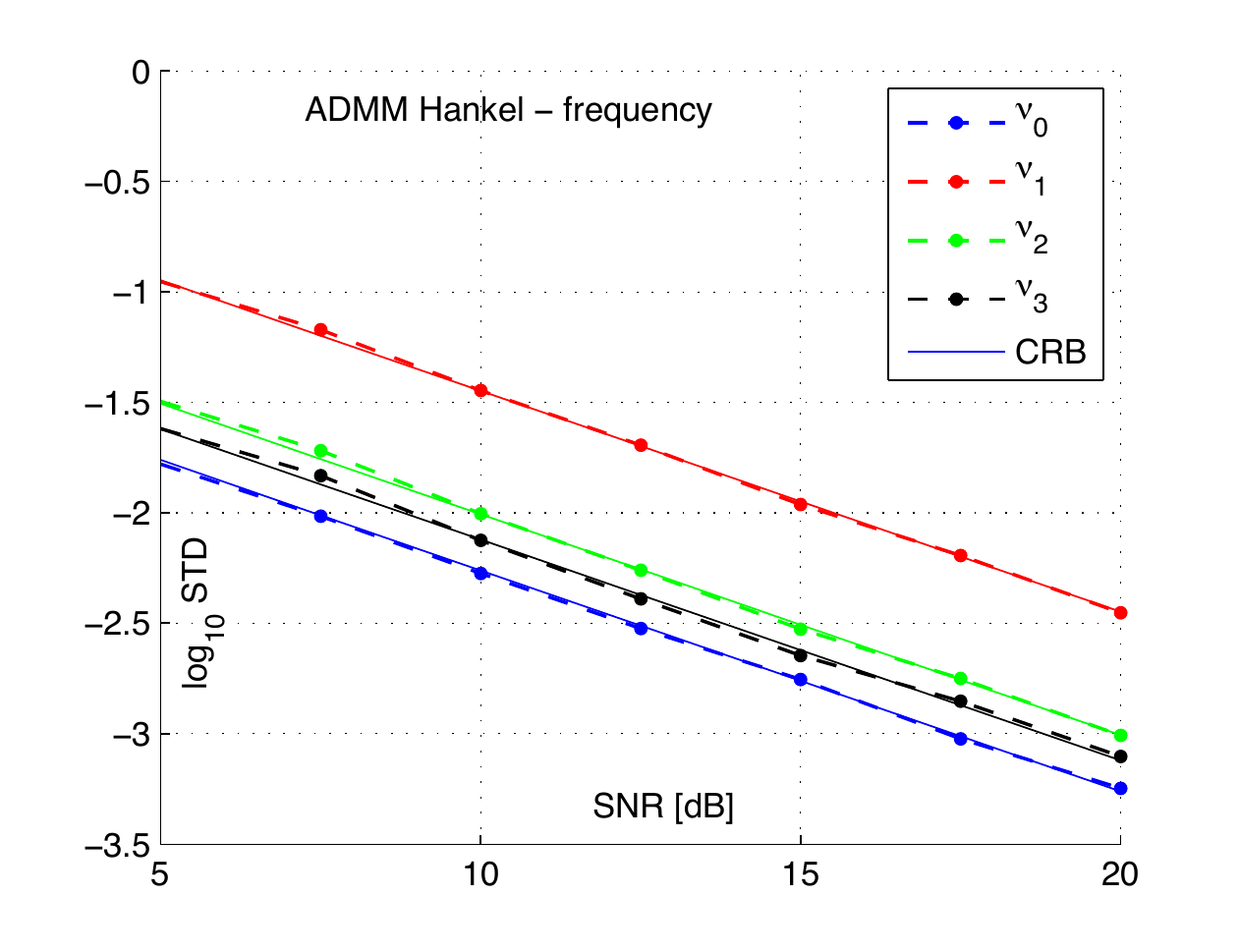}\includegraphics[width=0.5\linewidth]{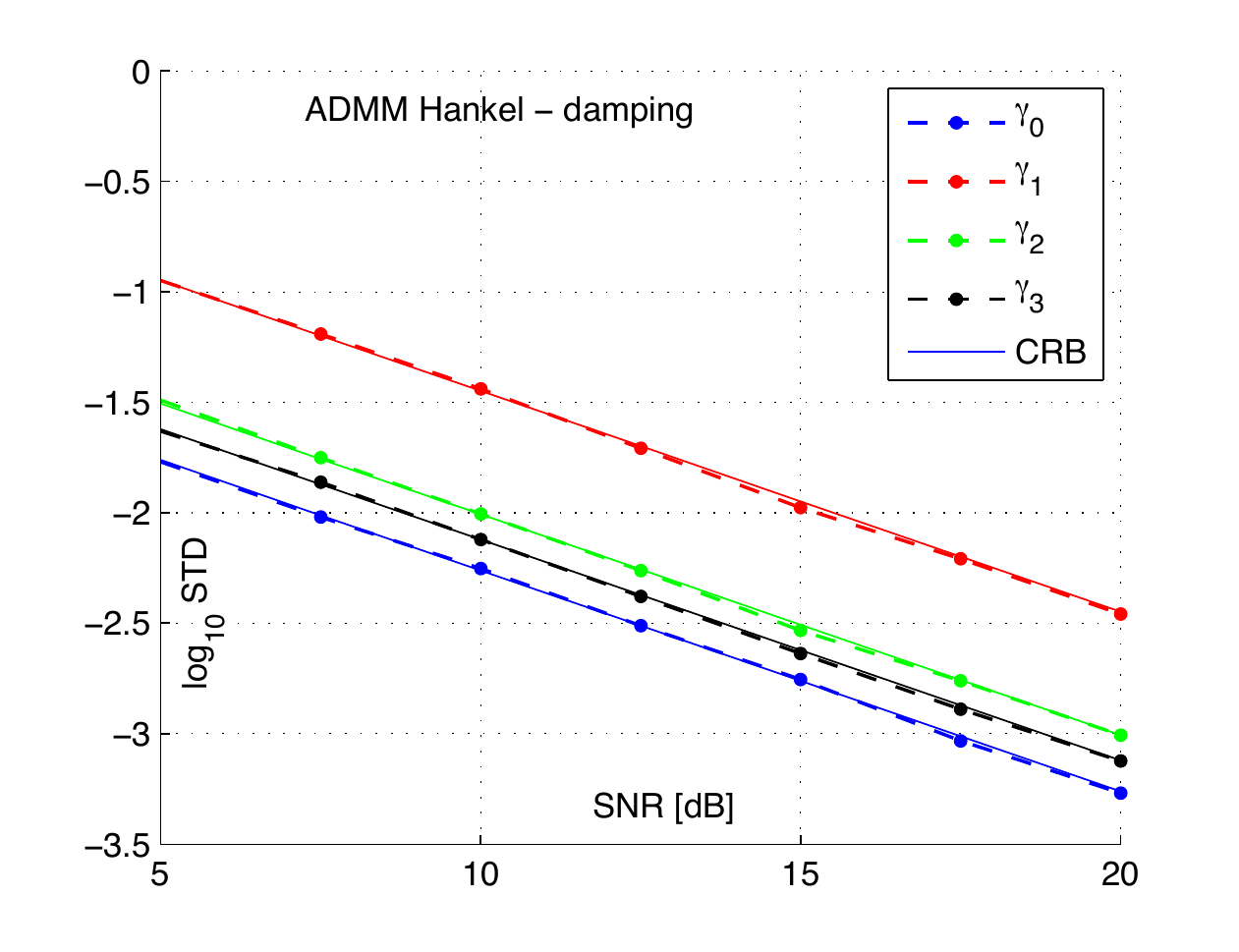}
\caption{\label{fig:perf} Performance for the estimation of $\freq$ (left) and $\damp$ (right), $p=0,\cdots,3$, for  ESPRIT (top) and the ADMM procedure (bottom). The CRBs are plotted in solid lines.}
\end{figure}
\begin{figure}[tb]
\centering
\includegraphics[width=0.5\linewidth]{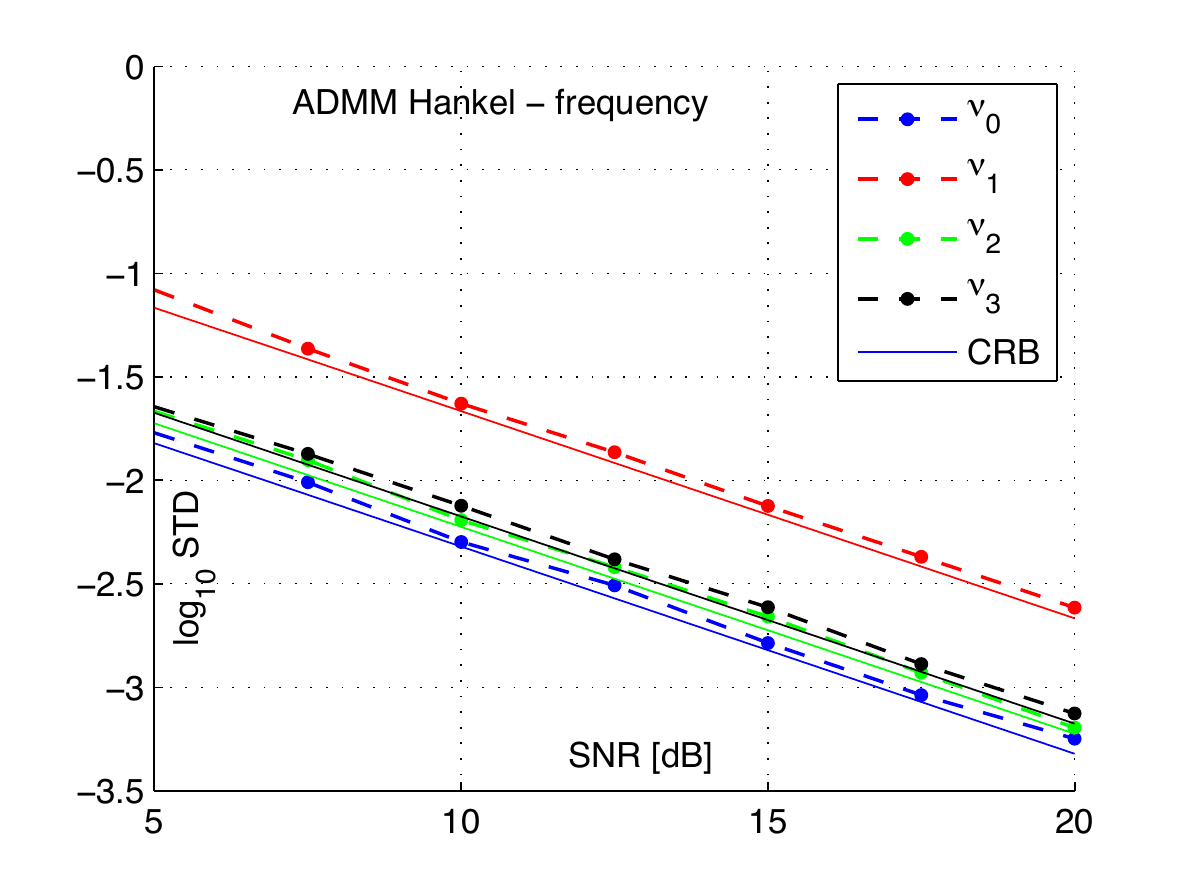}\includegraphics[width=0.5\linewidth]{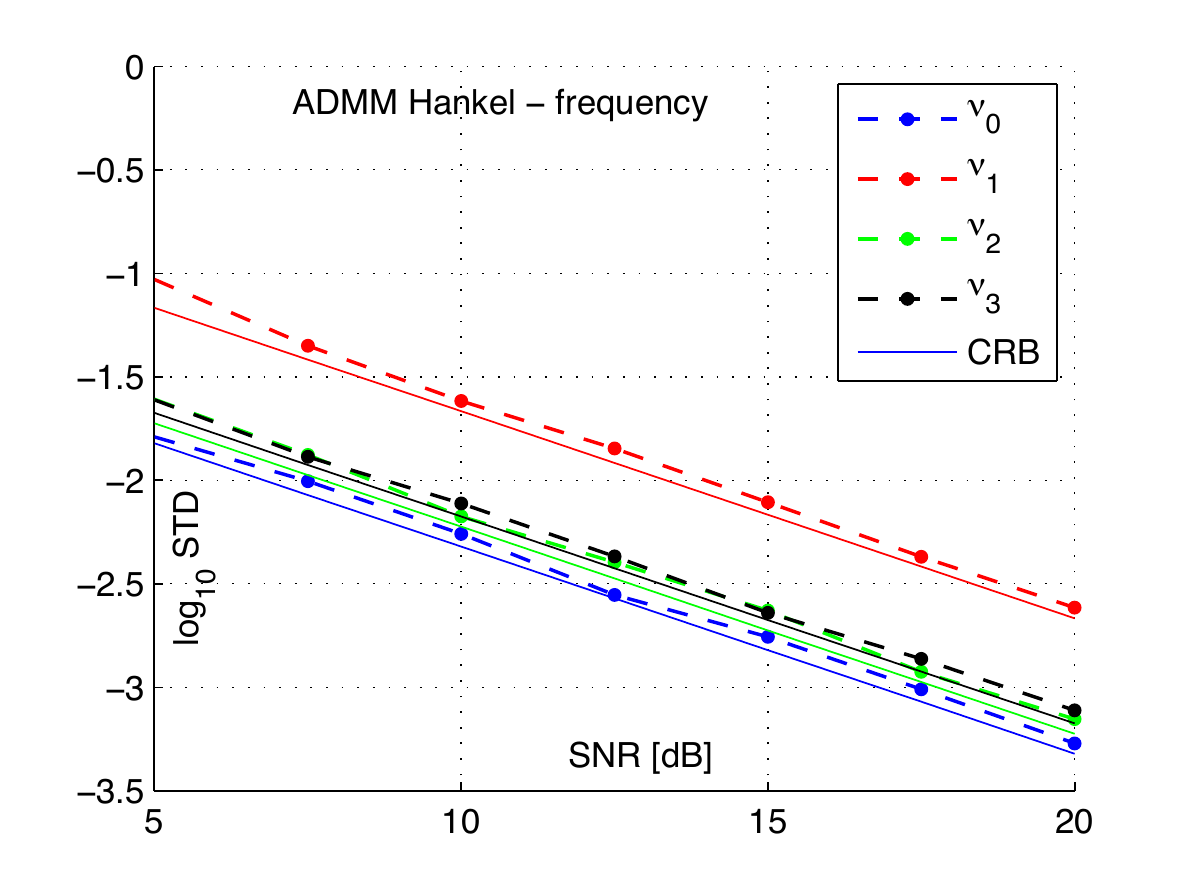}
\caption{\label{fig:perfmissing} Performance for the estimation of $\freq$, $p=0,\cdots,3$, for the ADMM procedure when $64$ out of $321$ samples are missing at random positions (left) and as one consecutive block (right). The solid lines indicate the CRBs.}
\end{figure}
\begin{figure}[tb]
\includegraphics[width=\linewidth]{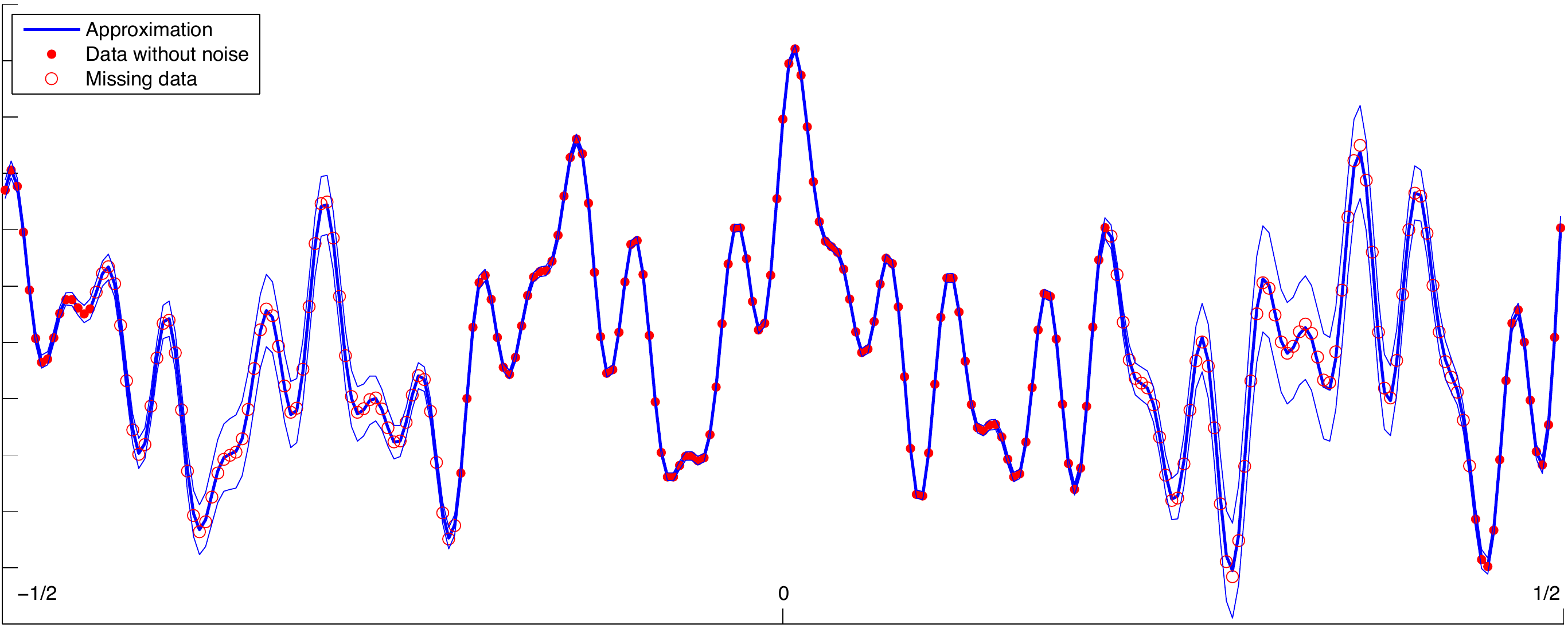}
\caption{\label{fig:appro}Approximation of $P=7$ exponentials for $2N+1=257$ samples of which two blocks of $60$ samples are missing ($SNR=10dB$). Real parts of noise-free available samples (red dots) and missing samples (red circles); mean and $\pm 1.96$ standard deviation error band (blue solid lines).
}
\end{figure}

\section{Estimation performance}
\label{sec:results}

We analyze the estimation performance of the proposed frequency estimation algorithm by considering different numerical simulations conducted with $P=4$ and $P=7$ complex exponentials embedded in circular white Gaussian noise ($500$ independent realisations) for different signal-to-noise ratios (SNR). The model parameters $\boldsymbol{c}$ and $\boldsymbol{\zeta}$ are summarized in Tab. \ref{tab:param} (units correspond to uniform sampling of the interval $t\in[-\frac{1}{2},\frac{1}{2}]$). The performance are compared with the theoretical Cram\'er-Rao bounds (CRBs) for the estimation problem (see, e.g., \cite{Gudmundson2012,Stoica2005}).

\subsection{Estimation performance for uniformly sampled data}
In Figure \ref{fig:perf}, the standard deviations (STDs) of the estimates of $\damp$ and $\freq$, $p=0,\cdots,3$, obtained with the ADMM procedure and with the ESPRIT method, are plotted as a function of SNR together with the (square root of the) CRBs. Note that the components $p=1$ and $p=2$ are close in frequency but significantly differ in amplitude (cf., Tab. \ref{tab:param}). 
We observe that the proposed ADMM procedure provides estimates attaining the corresponding CRBs, indicating that the procedure realizes the minimum variance estimator. It consistently outperforms the ESPRIT method in terms of STD, and especially so for small SNR and for the two components that are close in frequency. 

\subsection{Estimation performance for missing data}
In Fig. \ref{fig:perfmissing}, we summarize  the performance of the ADMM procedure for estimating $\freq$, $p=0,\cdots,3$, for $2N+1=257+64$ samples of which $64$ samples are missing 
at random positions (left) and as one consecutive block centred at random position (right). 
The STDs are again found to reach the (square root of the) CRBs, demonstrating that the ADMM procedure is not affected by missing data. Similar results are obtained for $\damp$ and are not reported here for space limitation reasons.

\subsection{Approximation for missing data} Finally, Figure \ref{fig:appro} illustrates the approximation performance of the ADMM procedure for $2N+1=257$ samples of $P=7$ exponentials of which two blocks of 60 consecutive samples are missing ($SNR=10dB$, mean over 100 realizations). Despite the fact that $47\%$ of the samples are missing, the obtained approximations are notably consistent with the underlying noise-free function. The maximum observed differences $|\hat\nu_p-\freq|$ are found to be below $10^{-3}$ times the sampling frequency.

\section{Conclusions}
\label{sec:conclusions}
A high-resolution parametric frequency estimation procedure that is based on approximation with sums of $P$ complex exponentials was proposed. In the proposed algorithm, the Kronecker theorem was used to cast the approximation problem in terms of generating functions for Hankel matrices of rank $P$. The resulting optimization problem was addressed by an ADMM procedure. The Hankel matrices obtained from the ADMM procedure were then used to compute the parameters of the complex exponential model. This is in contrast to other methods, such as classical NLS or subspace methods. Although the optimization problem considered in this paper is non-convex (and the problems of local minima are inherited with non-convexity, as is the case for classical NLS) numerical simulations indicated that the ADMM procedure yields excellent practical performance.

The ADMM procedure can be applied to equally spaced samples, including situations with missing data, and does not rely on explicit noise model assumptions. Despite its versatility, the resulting algorithm is simple and easy to implement. The method has been presented here for complex-valued data. However, its real-valued counterpart can be obtained in a similar fashion as outlined in this work.

The procedure will be applied to frequency node estimation in nuclear quadrupole resonance applications. Future work includes the extension to multiple time series and to the estimation of directions of arrival.

\section{Acknowledgment}
The work was supported by the Swedish Foundation for International
Cooperation in Research and Higher Education, the Swedish
Research Council and the Craaford Foundation. Part of this work was conducted within the Labex CIMI during visits of F. Andersson at University of Toulouse.

\bibliographystyle{plain}

\begin{thebibliography}{10}

\bibitem{AAK}
V.~M. Adamjan, D.~Z. Arov, and M.~G. Kre{\u\i}n.
\newblock Infinite {H}ankel matrices and generalized problems of
  {C}arath\'eodory-{F}ej\'er and {F}. {R}iesz.
\newblock {\em Funkcional. Anal. i Prilo\v zen.}, 2(1):1--19, 1968.

\bibitem{AlternatingProjections}
F.~{Andersson} and M.~{Carlsson}.
\newblock {Alternating projections on non-tangential manifolds}.
\newblock {\em arXiv:1107.4055}, 2011.

\bibitem{JAT}
F.~Andersson, M.~Carlsson, and M.V. de~Hoop.
\newblock Sparse approximation of functions using sums of exponentials and aak
  theory.
\newblock {\em Journal of Approximation Theory}, 163(2):213--248, February
  2011.

\bibitem{Beylkin_Monzon_2005}
G.~Beylkin and L.~Monzon.
\newblock On approximation of functions by exponential sums.
\newblock {\em Applied and Computational Harmonic Analysis}, 19(1):17--48, July
  2005.

\bibitem{Boyd_ADMM}
S.~Boyd, N.~Parikh, E.~Chu, B.~Peleato, and J.~Eckstein.
\newblock Distributed optimization and statistical learning via the alternating
  direction method of multipliers.
\newblock {\em Foundations and Trends{\textregistered} in Machine Learning},
  3(1):1--122, 2011.

\bibitem{NLLS1986}
Y.~Bresler and A.~Macovski.
\newblock Exact maximum likelihood parameter estimation of superimposed
  exponential signals in noise.
\newblock {\em IEEE Trans. Acoustic, Speech and Signal Process.},
  34(5):1081--1089, 1986.

\bibitem{Gudmundson2012}
E.~Gudmundson, P.~Wirfalt, A.~Jakobsson, and M.~Jansson.
\newblock An esprit-based parameter estimator for spectroscopic data.
\newblock In {\em IEEE Statistical Signal Processing Workshop (SSP'12)}, pages
  77 -- 80, 2012.

\bibitem{Horn}
R.~A. Horn and C.~R. Johnson.
\newblock {\em Topics in matrix analysis}.
\newblock Cambridge University Press, Cambridge, 1994.

\bibitem{MinNorm1983}
R.~Kumaresan and D.~W. Tufts.
\newblock Estimating the angles of arrival of multiple plane waves.
\newblock {\em IEEE Trans. Aerospace and Electronic Systems}, 19:134--139,
  1983.

\bibitem{ESPRIT}
R.~Roy and T.~Kailath.
\newblock {ESPRIT-estimation of signal parameters via rotational invariance
  techniques}.
\newblock {\em IEEE Trans. on Acoustics, Speech and Signal Process.},
  37(7):984--995, 1989.

\bibitem{Schmidt}
R.~Schmidt.
\newblock Multiple emitter location and signal parameter estimation.
\newblock {\em IEEE Trans. on Antennas and Propagation}, 34(3):276 -- 280,
  1986.

\bibitem{Stoica2005}
P.~Stoica and R.~Moses.
\newblock {\em Spectral analysis of Signals}.
\newblock Prentice--Hall, 2005.

\end{thebibliography}

\end{document}